\documentclass[12pt]{article}
\usepackage{amsmath}
\usepackage{amssymb}
\usepackage{amsfonts}
\usepackage{epsfig}

\begin{document}

\title{Stochastic solutions of nonlinear pde's: McKean versus superprocesses}
\author{R. Vilela Mendes \\
CMAF - Complexo Interdisciplinar, Universidade de Lisboa\\
(Av. Gama Pinto 2, 1649-003, Lisbon, vilela@cii.fc.ul.pt)\\
Instituto de Plasmas e Fus\~{a}o Nuclear, IST\\
(Av. Rovisco Pais, 1049-001 Lisbon)}
\date{ \ }
\maketitle

\begin{abstract}
Stochastic solutions not only provide new rigorous results for nonlinear
pde's but also, through its local non-grid nature, are a natural tool for
parallel computation.

There are two methods to construct stochastic solutions: the McKean method
and superprocesses. Here a comparison is made of these two approaches and
their strenghts and limitations are discussed.
\end{abstract}

\section{Introduction: Stochastic solutions and their uses}

A stochastic solution of a linear or nonlinear partial differential equation
is a stochastic process which, when started from a particular point in the
domain generates after a time $t$ a boundary measure which, when integrated
over the initial condition at $t=0$, provides the solution at the point $x$
and time $t$. For example for the heat equation%
\begin{equation}
\partial _{t}u(t,x)=\frac{1}{2}\frac{\partial ^{2}}{\partial x^{2}}%
u(t,x)\qquad \text{\textnormal{with}}\qquad u(0,x)=f(x)  \label{1.1}
\end{equation}%
the stochastic process is Brownian motion and the solution is 
\begin{equation}
u(t,x)=\mathbb{E}_{x}f(X_{t})  \label{1.2}
\end{equation}%
$\mathbb{E}_{x}$ meaning the expectation value, starting from $x$, of the
process%
\begin{equation}
dX_{t}=dB_{t}  \label{1.3}
\end{equation}%
The domain here is $\mathbb{R}\times \left[ 0,t\right) $ and the expectation
value in (\ref{1.2}) is indeed the inner product $\left\langle \mu
_{t},f\right\rangle $ of the initial condition $f$ with the measure $\mu
_{t} $ generated by the Brownian motion at the $t-$boundary. The usual
integral solution,%
\begin{equation}
u\left( t,x\right) =\frac{1}{2\sqrt{\pi }}\int \frac{1}{\sqrt{t}}\exp \left(
-\frac{\left( x-y\right) ^{2}}{4t}\right) f\left( y\right) dy  \label{1.4}
\end{equation}%
with the heat kernel, has exactly the same interpretation. Of course, an
important condition for the stochastic process (Brownian motion in this
case) to be considered \textit{the} solution of the equation is the fact
that the same process works for any initial condition. This should be
contrasted with stochastic processes constructed from particular solutions.

That the solutions of linear elliptic and parabolic equations, both with
Cauchy and Dirichlet boundary conditions, have a probabilistic
interpretation is a classical result and a standard tool in potential theory 
\cite{Getoor} \cite{Bass1} \cite{Bass2}. In contrast with the linear
problems, explicit solutions in terms of elementary functions or integrals
for nonlinear partial differential equations are only known in very
particular cases. Therefore the construction of solutions through stochastic
processes, for nonlinear equations, has become an active field in recent
years. The first stochastic solution for a nonlinear pde was constructed by
McKean \cite{McKean} for the KPP equation. Later on, the exit measures
provided by diffusion plus branching processes \cite{Dynkin1} \cite{Dynkin2}
as well as the stochastic representations recently constructed for the
Navier-Stokes \cite{Jan} \cite{Waymire} \cite{Bhatta1} \cite{Ossiander} \cite%
{Orum}, the Vlasov-Poisson \cite{Vilela1} \cite{Vilela2} \cite{Vilela4}, the
Euler \cite{Vilela3} and a fractional version of the KPP equation \cite%
{Cipriano} define solution-independent processes for which the mean values
of some functionals are solutions to these equations. Therefore, they are 
\textit{exact stochastic solutions}.

In the stochastic solutions one deals with a process that starts from the
point where the solution is to be found, a functional being then computed on
the boundary or in some cases along the whole sample path. In addition to
providing new exact results, the stochastic solutions are also a promising
tool for numerical implementation. This is because:

(i) Deterministic algorithms grow exponentially with the dimension $d$ of
the space , roughly $N^{d}$ ($\frac{L}{N}$ being the linear size of the
grid). This implies that to have reasonable computing times, the number of
grid points may not be sufficient to obtain a good local resolution for the
solution. In contrast a stochastic simulation only grows with the dimension
of the process, typically of order $d$.

(ii) In general, deterministic algorithms aim at obtaining the global
behavior of the solution in the whole domain. That means that even if an
efficient deterministic algorithm exists for the problem, a stochastic
algorithm might still be competitive if only localized values of the
solution are desired. This comes from the very nature of the stochastic
representation processes that always start from a definite point of the
domain. According to what is desired, real or Fourier space representations
should be used. For example by studying only a few high Fourier modes one
may obtain information on the small scale fluctuations that only a very fine
grid would provide in a deterministic algorithm.

(iii) Each time a sample path of the process is implemented, it is
independent from any other sample paths that are used to obtain the
expectation value. Likewise, paths starting from different points are
independent from each other. Therefore the stochastic algorithms are a
natural choice for parallel and distributed implementation. Provided some
differentiability conditions are satisfied, the process also handles equally
well simple or complex boundary conditions.

(iv) Stochastic algorithms may also be used for domain decomposition
purposes \cite{Acebron1} \cite{Acebron2} \cite{Acebron3}. One may, for
example, decompose the space in subdomains and then use in each one a
deterministic algorithm with Dirichlet boundary conditions, the values on
the boundaries being determined by a stochastic algorithm, thus minimizing
the time-consuming communication problem between domains.

There are basically two methods to construct stochastic solutions. The first
method, which will be called the McKean method, is essentially a
probabilistic interpretation of the Picard series. The differential equation
is written as an integral equation which is rearranged in a such a way that
the coefficients of the successive terms in the Picard iteration obey a
normalization condition. The Picard iteration is then interpreted as an
evolution and branching process, the stochastic solution being equivalent to
importance sampling of the normalized Picard series. The second method
constructs the boundary measures of a measure-valued stochastic process (a
superprocess) and obtain the solutions of the differential equation by a
scaling procedure.

In this paper the two methods are compared and in the final sections one
shows how the superprocess construction may be extended to larger classes of
partial differential equations by going from process on measures to
processes on signed measures and processes on distributions.

\section{McKean and superprocesses: The KPP\ equation}

\subsection{The KPP equation: McKean's formulation}

To illustrate the two methods for the construction of stochastic solutions a
classical example will be used, namely the KPP equation \cite{McKean}

\begin{equation}
\frac{\partial v}{\partial t}=\frac{1}{2}\frac{\partial ^{2}v}{\partial x^{2}%
}+v^{2}-v  \label{2.1}
\end{equation}%
with initial data $v\left( 0,x\right) =g\left( x\right) $

Let $G\left( t,x\right) $ be the Green's operator for the heat equation $%
\partial _{t}v(t,x)=\frac{1}{2}\frac{\partial ^{2}}{\partial x^{2}}v(t,x)$%
\begin{equation}
G\left( t,x\right) =e^{\frac{1}{2}t\frac{\partial ^{2}}{\partial x^{2}}}
\label{2.2}
\end{equation}%
Then the equation in integral form is%
\begin{equation}
v\left( t,x\right) =e^{-t}G\left( t,x\right) g\left( x\right)
+\int_{0}^{t}e^{-\left( t-s\right) }G\left( t-s,x\right) v^{2}\left(
s,x\right) ds  \label{2.3}
\end{equation}%
Denoting by $\left( \xi _{t},\Pi _{x}\right) $ a Brownian motion started
from time zero and coordinate $x$, Eq.(\ref{2.3}) may be rewritten%
\begin{eqnarray}
v\left( t,x\right) &=&\Pi _{x}\left\{ e^{-t}g\left( \xi _{t}\right)
+\int_{0}^{t}e^{-\left( t-s\right) }v^{2}\left( s,\xi _{t-s}\right)
ds\right\}  \notag \\
&=&\Pi _{x}\left\{ e^{-t}g\left( \xi _{t}\right)
+\int_{0}^{t}e^{-s}v^{2}\left( t-s,\xi _{s}\right) ds\right\}  \label{2.3a}
\end{eqnarray}%
Therefore the solution is obtained by the following process:

At the initial time, a single particle begins a Brownian motion, starting
from $x$ and continuing for an exponential holding time $T$ with $P\left(
T>t\right) =e^{-t}$.

Then, at $T$, the particle splits into two, the new particles continuing
along independent Brownian paths starting from $x\left( T\right) $. These
particles, in turn, are subjected to the same splitting rule, meaning that
after an elapsed time $t>0$ one has $n$ particles located at $x_{1}\left(
t\right) ,x_{2}\left( t\right) ,\cdots x_{n}\left( t\right) $ with $P\left(
n=k\right) =e^{-t}\left( 1-e^{-t}\right) ^{k-1}$.

The solution of (\ref{2.3}) is obtained by%
\begin{equation}
v\left( t,x\right) =\mathbb{E}\left\{ g\left( x_{1}(t)\right) g\left(
x_{2}(t)\right) \cdots g\left( x_{n}(t)\right) \right\}  \label{2.4}
\end{equation}

\begin{figure}[htb]
\begin{center}
\psfig{figure=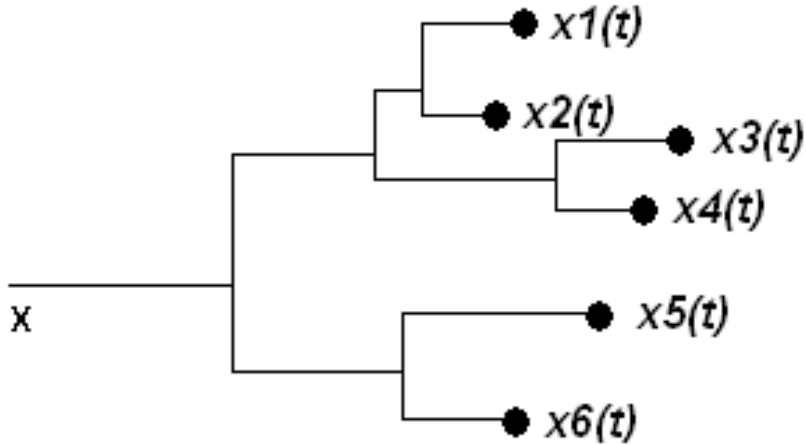,width=11truecm}
\end{center}
\caption{The McKean process for the KPP equation}
\end{figure}

An equivalent interpretation consists in considering the process propagating
backwards in time from time $t$ at the point $x$ and, when it reaches time
zero, it samples the initial condition. That is, the process generates a
measure at the $t=0$ boundary which is then applied to the function $g\left(
x\right) =v\left( 0,x\right) $.

This construction, which expresses the solution as a stochastic
multiplicative functional of the initial condition, is also qualitatively
equivalent to importance sampling of the Picard iteration of Eq.(\ref{2.3}).

A sufficient condition for the existence of (\ref{2.4}) is $\left\vert
g\left( x\right) \right\vert \leq 1$ or, almost surely, $\left\vert g\left(
x\right) \right\vert \left( 1-e^{-t}\right) \leq 1$.

Another probabilistic approach to this type of equations is through the
construction of superprocesses. In many cases a superprocess may be looked
at as the scaling limit of a branching particle system. The point of view
used in the derivation of superprocesses is different from the derivation
above. In the next subsection a short introduction to superprocesses is
sketched and then the KPP\ solution is constructed via superprocesses.

\subsection{Branching exit measures and superprocesses}

Let $\left( E,\mathcal{B}\right) $ be a measurable space and $M_{+}\left(
E\right) $ the space of finite measures in $E$. Denote by $\left(
X_{t},P_{0\,,\mu }\right) $ a branching stochastic process with values in $%
M_{+}\left( E\right) $ and transition probability $P_{0,\mu }$ starting from
time $0$ and measure $\mu $. The process is said to satisfy the \textit{%
branching property} if given $\mu =\mu _{1}+\mu _{2}$%
\begin{equation}
P_{0,\mu }=P_{0,\mu _{1}}\ast P_{0,\mu _{2}}  \label{2.5}
\end{equation}%
that is, after the branching $\left( X_{t}^{1},P_{0,\mu _{1}}\right) $ and $%
\left( X_{t}^{2},P_{0\,,\mu _{2}}\right) $ are independent and $%
X_{t}^{1}+X_{t}^{2}$ has the same law as $\left( X_{t},P_{r0,\mu }\right) $.
In terms of the \textit{transition operator }$V_{t}$ operating on functions
on $E$ this is%
\begin{equation}
V_{t}f\left( \mu _{1}+\mu _{2}\right) =V_{t}f\left( \mu _{1}\right)
+V_{t}f\left( \mu _{2}\right)  \label{2.6}
\end{equation}%
where $e^{-\left\langle V_{t}f,\mu \right\rangle }\circeq P_{r,\mu
}e^{-\left\langle f,X_{t}\right\rangle }$ or%
\begin{equation}
V_{t}f\left( \mu \right) =-\log P_{0,\mu }e^{-\left\langle
f,X_{t}\right\rangle }  \label{2.7}
\end{equation}%
$V_{t}$ is called the \textit{log-Laplace semigroup} associated to $X_{t}$.
In (\ref{2.7}) if the initial measure $\mu $ is $\delta _{x}$ one writes%
\begin{equation}
V_{t}f\left( x\right) =-\log P_{0,x}e^{-\left\langle f,X_{t}\right\rangle }
\label{2.8}
\end{equation}%
By (\ref{2.5}) the probability law of $X_{t}$ is infinitely divisible.

Now in $S=\left[ 0,\infty \right) \times E$ consider a set $Q\subset S$ and
the associated branching exit process \ $\left( X_{Q},P_{\mu }\right) $
composed of a propagating Markov process in $E$, $\xi =\left( \xi _{t},\Pi
_{0,x}\right) $, a set of probabilities $p_{n}(t,x)$ describing the
branching and a parameter $k$ defining the lifetime.%
\begin{equation}
u\left( x\right) =V_{Q}f\left( x\right) =-\log P_{0,x}e^{-\left\langle
f,X_{Q}\right\rangle }  \label{2.10}
\end{equation}%
$\left\langle f,X_{Q}\right\rangle $ is the integral of the function $f$ on
the (space-time) boundary with the boundary exit measure generated by the
process. One says that this branching exit process is a $\left( \xi ,\psi
\right) -$\textit{superprocess} if $u\left( x\right) $ satisfies the equation%
\begin{equation}
u+G_{Q}\psi \left( u\right) =K_{Q}f  \label{2.11}
\end{equation}%
where $G_{Q}$ is the Green operator,%
\begin{equation}
G_{Q}f\left( r,x\right) =\Pi _{0,x}\int_{0}^{\tau }f\left( s,\xi _{s}\right)
ds  \label{2.12}
\end{equation}%
$K_{Q}$ the Poisson operator%
\begin{equation}
K_{Q}f\left( x\right) =\Pi _{0,x}1_{\tau <\infty }f\left( \xi _{\tau }\right)
\label{2.13}
\end{equation}%
$\psi \left( u\right) $ means $\psi \left( 0,x;u\left( 0,x\right) \right) $
and $\tau $ is the exit time from $Q$.

The superprocess is constructed as follows: Let $\varphi \left( s,x;z\right) 
$ be the offspring generating function at time $s$ and point $x$%
\begin{equation}
\varphi \left( s,x;z\right) =c\sum_{0}^{\infty }p_{n}(s,x)z^{n}  \label{2.14}
\end{equation}%
where $\sum_{n}p_{n}=1$ and $c$ denotes the branching intensity.

Then for $e^{-w\left( 0,x\right) }\circeq P_{0,x}e^{-\left\langle
f,X_{Q}\right\rangle }$ one has%
\begin{equation}
P_{0,x}e^{-\left\langle f,X_{Q}\right\rangle }\circeq e^{-w\left( 0,x\right)
}=\Pi _{0,x}\left[ e^{-k\tau }e^{-f\left( \tau ,\xi _{\tau }\right)
}+\int_{0}^{\tau }dske^{-ks}\varphi \left( s,\xi _{s};e^{-w\left( \tau
-s,\xi _{s}\right) }\right) \right]  \label{2.15}
\end{equation}%
The measure-valued process starts from $\delta _{x}$ at time $0$, $\tau $ is
the first exit time from $Q$ and $f\left( \tau ,\xi _{\tau }\right) $ the
value of a function in the boundary $\partial Q$.

Using $\int_{0}^{\tau }ke^{-ks}ds=1-e^{-k\tau }$ and the Markov property $%
\Pi _{0,x}1_{s<\tau }\Pi _{s,\xi _{s}}=\Pi _{0,x}1_{s<\tau }$, Eq.(\ref{2.15}%
) for $e^{-w\left( 0,x\right) }$ is converted into%
\begin{equation}
e^{-w\left( 0,x\right) }=\Pi _{0,x}\left[ e^{-f\left( \tau ,\xi _{\tau
}\right) }+k\int_{0}^{\tau }ds\left[ \varphi \left( s,\xi _{s};e^{-w\left(
\tau -s,\xi _{s}\right) }\right) -e^{-w\left( \tau -s,\xi _{s}\right) }%
\right] \right]  \label{2.16}
\end{equation}%
This is lemma 1.2 in ch.4 of Ref.\cite{Dynkin1}. Because of the central role
of this result for the construction of superprocesses, a proof is included
in the Appendix with the notations used in this paper.

Eq.(2.11) is now obtained by a limiting process. Let in (\ref{2.16}) replace 
$w\left( 0,x\right) $ by $\beta w_{\beta }\left( 0,x\right) $ and $f$ by $%
\beta f$. $\beta $ is interpreted as the mass of the particles and when the
measure-valued process $X_{Q}\rightarrow \beta X_{Q}$ then $P_{\mu
}\rightarrow P_{\frac{\mu }{\beta }}$.%
\begin{equation}
e^{-\beta w\left( 0,x\right) }=\Pi _{0,x}\left[ e^{-\beta f\left( \tau ,\xi
_{\tau }\right) }+k_{\beta }\int_{0}^{\tau }ds\left[ \varphi _{\beta }\left(
s,\xi _{s};e^{-\beta w\left( \tau -s,\xi _{s}\right) }\right) -e^{-\beta
w\left( \tau -s,\xi _{s}\right) }\right] \right]  \label{2.17}
\end{equation}%
Defining%
\begin{equation}
u_{\beta }=\left( 1-e^{-\beta w_{\beta }}\right) /\beta \hspace{0.3cm};%
\hspace{0.3cm}f_{\beta }=\left( 1-e^{-\beta f}\right) /\beta  \label{2.18}
\end{equation}%
and%
\begin{equation}
\psi _{\beta }\left( r,x;u_{\beta }\right) =\frac{k_{\beta }}{\beta }\left(
\varphi \left( r,x;1-\beta u_{\beta }\right) -1+\beta u_{\beta }\right)
\label{2.19}
\end{equation}%
one obtains from (\ref{2.17})%
\begin{equation}
u_{\beta }\left( 0,x\right) +\Pi _{0,x}\int_{0}^{\tau }ds\psi _{\beta
}\left( s,\xi _{s};u_{\beta }\right) =\Pi _{0,x}f_{\beta }\left( \tau ,\xi
_{\tau }\right)  \label{2.20}
\end{equation}%
that is%
\begin{equation}
u_{\beta }+G_{Q}\psi _{\beta }\left( u_{\beta }\right) =K_{Q}f_{\beta }
\label{2.21}
\end{equation}%
When $\beta \rightarrow 0$, $f\rightarrow f_{\beta }$ and if $\psi _{\beta }$
goes to a well defined limit $\psi $ then $u_{\beta }$ tends to a limit $u$
solution of (\ref{2.11}) associated to a superprocess. Also one sees from (%
\ref{2.18}) that in the $\beta \rightarrow 0$ limit%
\begin{equation}
u_{\beta }\rightarrow w_{\beta }=-\log P_{0,x}e^{-\left\langle
f,X_{Q}\right\rangle }  \label{2.21a}
\end{equation}%
as in Eq.(\ref{2.10}). The superprocess corresponds to a cloud of particles
for which both the mass and the lifetime tend to zero.

\subsection{The KPP\ equation as a superprocess}

When the integral Eq.(\ref{2.3}) is interpreted probabilistically, it may be
identified with Eq.(\ref{2.15}) with $k=1$, $e^{-w\left( 0,x\right)
}=v\left( \tau ,x\right) $, $e^{-f\left( \tau ,\xi _{\tau }\right) }=g\left(
\xi _{\tau }\right) $, $\varphi \left( s,\xi _{s};e^{-w\left( \tau -s,\xi
_{s}\right) }\right) =v^{2}\left( \tau -s,\xi _{s}\right) $. Therefore the
McKean probabilistic construction corresponds to an intermediate step in the
superprocess construction. At this level the process that is considered in
Eq.(\ref{2.15}) is the same as in McKean's construction. Summing over the
exit measure, the solution is%
\begin{equation}
v\left( t,x\right) =e^{-\left\langle f,X_{Q}\right\rangle
}=e^{-\sum_{i}f\left( \xi _{\tau _{i}}\right) }=e^{\sum_{i}\log g\left( \xi
_{\tau _{i}}\right) }=\Pi _{i}g\left( \xi _{\tau _{i}}\right)  \label{2.22}
\end{equation}%
essentially the same as in (\ref{2.4}). However, there are two differences.
First, the initial condition $g$ must be positive to have a well-defined
logarithm. This is a restriction as compared to McKean's construction. But,
on the other hand, the interpretation as an exit measure, allows to deal
with Cauchy problems with boundary conditions. The exit measure is from the
set $Q=\left[ 0,t\right] \times \Omega $, $\tau $ being the time at which
the $\left( t,\xi _{t}\right) $ process reaches $\partial \Omega $ or $\tau
=t$ inside $\Omega $.

For the superprocess, let $u\left( t,x\right) =1-v\left( t,x\right) $, which
satisfies the equation%
\begin{equation}
\frac{\partial u}{\partial t}=\frac{1}{2}\frac{\partial ^{2}u}{\partial x^{2}%
}-u^{2}+u  \label{2.23}
\end{equation}%
and the integral equation%
\begin{equation}
u\left( t,x\right) =G\left( t,x\right) \left( 1-g\left( x\right) \right)
+\int_{0}^{t}G\left( s,x\right) \left( u\left( t-s,x\right) -u^{2}\left(
t-s,x\right) \right) ds  \label{2.24}
\end{equation}%
or%
\begin{equation}
u\left( t,x\right) +\Pi _{x}\int_{0}^{t}\left( u^{2}\left( t-s,\xi
_{s}\right) -u\left( t-s,\xi _{s}\right) \right) ds=\Pi _{x}\left( 1-g\left(
\xi _{t}\right) \right)  \label{2.25}
\end{equation}%
that is for KPP%
\begin{equation}
\psi \left( 0,x;u\right) =u^{2}-u  \label{2.26}
\end{equation}%
Equating with (\ref{2.19}) one obtains%
\begin{eqnarray}
\psi _{\beta }\left( 0,x;u_{\beta }\right) &=&\frac{k_{\beta }}{\beta }%
\left( \varphi \left( 0,x;1-\beta u_{\beta }\right) -1+\beta u_{\beta
}\right)  \notag \\
&=&\frac{k_{\beta }}{\beta }\left( c\sum p_{n}\left( 1-\beta u_{\beta
}\right) ^{n}-1+\beta u_{\beta }\right)  \notag \\
&=&\frac{k_{\beta }c}{\beta }\left( \beta ^{2}u_{\beta }^{2}-\beta u_{\beta
}\right)  \notag \\
&=&u^{2}-u  \label{2.27}
\end{eqnarray}%
with $p_{n}=\delta _{n,2}$. Therefore $c=\beta =1$ and $k_{\beta }=1$. That
is, for KPP the superprocess is not a scaling limit. It coincides with the
McKean process. However in this case, because $\beta =1$ instead of $\beta
\rightarrow 0$, the solution is given by $\left( 1-e^{-w}\right) $ instead
of (\ref{2.10}).

However, the power of the superprocesses is that, with other limiting
choices of $\beta $, stochastic solutions may be constructed for other
equations, in particular for solutions without the natural Poisson clock
provided by the term $-v$ which is present in the KPP equation. For example
for%
\begin{equation}
\frac{\partial v}{\partial t}=\frac{1}{2}\frac{\partial ^{2}v}{\partial x^{2}%
}+v^{2}  \label{2.28}
\end{equation}%
with $u=-v$%
\begin{equation}
\frac{\partial u}{\partial t}=\frac{1}{2}\frac{\partial ^{2}u}{\partial x^{2}%
}-u^{2}  \label{2.28a}
\end{equation}%
one has%
\begin{equation}
\psi \left( 0,x;u\right) =u^{2}  \label{2.29}
\end{equation}%
which equated with (\ref{2.19})%
\begin{eqnarray}
\psi _{\beta }\left( 0,x;u_{\beta }\right) &=&\frac{k_{\beta }}{\beta }%
\left( \beta u_{\beta }-1+\sum_{n=0}^{2}p_{0}+p_{1}\left( 1-\beta u_{\beta
}\right) +p_{2}\left( 1-\beta u_{\beta }\right) ^{2}\right)  \notag \\
&=&u_{\beta }^{2}  \label{2.30}
\end{eqnarray}%
leads to%
\begin{equation}
p_{1}=0;\hspace{0.3cm}p_{0}=p_{2}=\frac{1}{2};\hspace{0.3cm}k_{\beta }=\frac{%
2}{\beta }  \label{2.31}
\end{equation}%
In this case one may let $\beta \rightarrow 0$, the solution is given by (%
\ref{2.10}) and the superprocess corresponds to the scaling limit $\left(
n\rightarrow \infty \text{ in Fig.2}\right) $ of a process where both the
mass and the lifetime of the particles tends to zero and at each bifurcation
point one has equal probability of either dying without offspring or having
two children (Fig.2)

\begin{figure}[htb]
\begin{center}
\psfig{figure=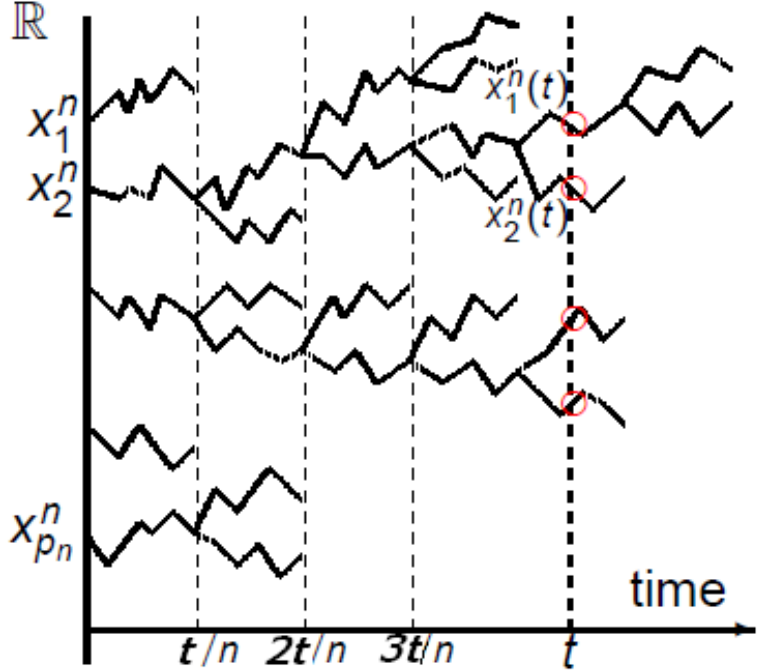,width=11truecm}
\end{center}
\caption{The branching process which in the scaling limit $n\rightarrow
\infty $ leads to the superprocess solution of Eq.(\protect\ref{2.28a})}
\end{figure}

This construction may be generalized for interactions $u^{\alpha }$ with $%
1<\alpha \leq 2$. With $z=1-\beta u_{\beta }$ one has%
\begin{eqnarray}
\varphi \left( 0,x;z\right) &=&\sum_{n}p_{n}z^{n}=z+\frac{\beta }{k_{\beta }}%
u_{\beta }^{\alpha }=z+\frac{\beta }{k_{\beta }}\frac{\left( 1-z\right)
^{\alpha }}{\beta ^{\alpha }}  \notag \\
&=&z+\frac{1}{k_{\beta }\beta ^{\alpha -1}}\left( 1-\alpha z+\frac{\alpha
\left( \alpha -1\right) }{2}z^{2}-\frac{\alpha \left( \alpha -1\right)
\left( \alpha -2\right) }{3!}z^{3}+\cdots \right)  \label{2.32}
\end{eqnarray}%
Choosing $k_{\beta }=\frac{\alpha }{\beta ^{\alpha -1}}$ the terms in $z$
cancel and for $1<\alpha \leq 2$ the coefficients of all the remaining $z$
powers are positive and may be interpreted as branching probabilities. It
would not be so for $\alpha >2$. Then%
\begin{equation}
p_{0}=\frac{1}{\alpha };\hspace{0.3cm}p_{1}=0;\hspace{0.3cm}\cdots \hspace{%
0.3cm}p_{n}=\frac{\left( -1\right) ^{n}}{\alpha }\left( 
\begin{array}{l}
\alpha \\ 
n%
\end{array}%
\right) \hspace{0.3cm}n\geq 2  \label{2.33}
\end{equation}%
with $\sum_{n}p_{n}=1$. With this choice of branching probabilities, $%
k_{\beta }=\frac{\alpha }{\beta ^{\alpha -1}}$ and $\beta \rightarrow 0$ one
obtains a superprocess which, through (\ref{2.10}), provides a solution to
the equation%
\begin{equation}
\frac{\partial u}{\partial t}=\frac{1}{2}\frac{\partial ^{2}u}{\partial x^{2}%
}-u^{\alpha }  \label{2.34}
\end{equation}%
for $1<\alpha \leq 2$.

However, the superprocess cannot be constructed for $\alpha >2$ because some
of the $z^{n}$ coefficients in the offspring generating function $\varphi
\left( 0,x;z\right) $ would be negative.

\section{Appendix: Proof of the lemma}

Let%
\begin{equation}
u\left( x,t\right) =\Pi _{0,x}\left\{ e^{-kt}u\left( \xi _{t},0\right)
+\int_{0}^{t}ke^{-ks}\Phi \left( \xi _{s},t-s\right) ds\right\}  \label{A.2}
\end{equation}

Then%
\begin{eqnarray}
\Pi _{0,x}\int_{0}^{t}ku\left( \xi _{s},t-s\right) ds &=&\Pi _{0,x}\left\{
\int_{0}^{t}ke^{-k\left( t-s\right) }u\left( \xi _{s+t-s},0\right) ds\right.
\notag \\
&&\left. +\int_{0}^{t}kds\int_{0}^{t-s}kds^{\prime }e^{-ks^{\prime }}\Phi
\left( \xi _{s+s^{\prime }},t-s-s^{\prime }\right) \right\}  \notag \\
&&  \label{A.3}
\end{eqnarray}%
Summing (\ref{A.2}) and (\ref{A.3})%
\begin{eqnarray}
&&u\left( x,t\right) +\Pi _{0,x}\int_{0}^{t}ku\left( \xi _{s},t-s\right) ds 
\notag \\
&=&\Pi _{0,x}\left\{ \left( e^{-kt}+\int_{0}^{t}ke^{-k\left( t-s\right)
}ds\right) u\left( \xi _{t},0\right) \right.  \notag \\
&&\left. +k\int_{0}^{t}e^{-ks}\Phi \left( \xi _{s},t-s\right)
ds+k\int_{0}^{t}ds\int_{0}^{t-s}kds^{\prime }e^{-ks^{\prime }}\Phi \left(
\xi _{s+s^{\prime }},t-s-s^{\prime }\right) ds^{\prime }\right\}  \notag \\
&&  \label{A.4}
\end{eqnarray}%
Changing variables in the last integral in (\ref{A.4}) from $\left(
s,s^{\prime }\right) $ to $\left( s,\sigma =s+s^{\prime }\right) $ one
obtains for the last term%
\begin{equation*}
k\int_{0}^{t}d\sigma \int_{0}^{\sigma }kdse^{-k\left( \sigma -s\right) }\Phi
\left( \xi _{\sigma },t-\sigma \right) ds
\end{equation*}%
and finally%
\begin{eqnarray}
&&u\left( x,t\right) +\Pi _{0,x}k\int_{0}^{t}u\left( \xi _{s},t-s\right) ds 
\notag \\
&=&\Pi _{0,x}\left\{ u\left( \xi _{t},0\right) +k\int_{0}^{t}\Phi \left( \xi
_{s},t-s\right) ds\right\}  \label{A.5}
\end{eqnarray}


\begin{thebibliography}{99}
\bibitem{Getoor} R. M. Blumenthal and R. K. Getoor; \textit{Markov processes
and potential theory}, Academic Press, New York 1968.

\bibitem{Bass1} R. F. Bass; \textit{Probabilistic techniques in analysis},
Springer, New York 1995.

\bibitem{Bass2} R. F. Bass; \textit{Diffusions and elliptic operators},
Springer, New York 1998.

\bibitem{McKean} H. P. McKean; Comm. on Pure and Appl. Math. 28 (1975)
323-331, 29 (1976) 553-554.

\bibitem{Dynkin1} E. B. Dynkin; \textit{Diffusions, Superdiffusions and
Partial Differential Equations, }AMS Colloquium Pubs., Providence 2002.

\bibitem{Dynkin2} E. B.Dynkin; \textit{Superdiffusions and positive
solutions of nonlinear partial differential equations}, AMS ,
Providence.2004.

\bibitem{Jan} Y. LeJan and A. S. Sznitman ; Prob. Theory and Relat. Fields
109 (1997) 343-366.

\bibitem{Waymire} E. C. Waymire; Prob. Surveys 2 (2005) 1-32.

\bibitem{Bhatta1} R. N. Bhattacharya et al. ; Trans. Amer. Math. Soc. 355
(2003) 5003-5040

\bibitem{Ossiander} M. Ossiander ; Prob. Theory and Relat. Fields 133 (2005)
267-298.

\bibitem{Orum} J. C. Orum; \textit{Stochastic cascades and 2D Fourier
Navier-Stokes equations}, in \textit{Lectures on multiscale and
multiplicative processes,} www.maphysto.dk/publications/MPS-LN/2002/11.pdf

\bibitem{Vilela1} R. Vilela Mendes and F. Cipriano; Commun. Nonlinear
Science and Num. Simul. 13 (2008) 221-226 and 1736.

\bibitem{Vilela2} E. Floriani, R. Lima and R. Vilela Mendes; European
Physical Journal D 46 (2008) 295-302 and 407.

\bibitem{Vilela3} R. Vilela Mendes; Stochastics 81 (2009) 279-297.

\bibitem{Vilela4} R. Vilela Mendes; J. Math. Phys. 51 (2010) 043101.

\bibitem{Cipriano} F. Cipriano, H. Ouerdiane and R. Vilela Mendes; Fract.
Calc. Appl. Anal.\textbf{\ }12 (2009) 47-56.

\bibitem{Acebron1} J. A. Acebr\'{o}n, A. Rodriguez-Rozas and R. Spigler; J.
of Comp. Physics 228 (2009) 5574--5591.

\bibitem{Acebron2} J.A. Acebr\'{o}n, A. Rodr\'{\i}guez-Rozas and R. Spigler;
J. on Scientific Computing 43 (2010) 135-157.

\bibitem{Acebron3} J.A. Acebr\'{o}n and A. Rodr\'{\i}guez-Rozas; J. Comp.
Phys. 230 (2011) 7891-7909.
\end{thebibliography}
\end{document}